\documentclass[11pt]{article}

\usepackage[bottom]{footmisc}

\usepackage[latin1]{inputenc}
\usepackage[all]{xypic}

\usepackage{amsmath,amssymb,amsfonts,enumerate,bbm,eucal}

\newcommand{\tmop}[1]{\operatorname{#1}}

\newcommand{\toisom}{\stackrel{\sim}{\to}}

\newcommand{\tensor}{\otimes}

\newcommand{\ideal}[1]{{\mathfrak #1}}

\newcommand{\id}{{\tmop{id}}}

\newcommand{\depth}{\tmop{depth}}

\newenvironment{enumeratealpha}{\begin{enumerate}[a{\textup{)}}]}{\end{enumerate}}

\DeclareMathOperator{\Ext}{Ext}
\DeclareMathOperator{\Tor}{Tor}

\DeclareMathOperator{\Ass}{Ass}

\DeclareMathOperator{\gldh}{gldh}

\DeclareMathOperator{\pd}{pd}

\DeclareMathOperator{\gr}{gr}

\DeclareMathOperator{\tot}{tot}

\newtheorem{lemma}{Lemma}[section]
\newtheorem{definition}{Definition}[section]
\newtheorem{proposition}{Proposition}[section]

\author{J\"urgen B\"ohm\thanks{jboehm@gmx.net}}
\title{A new proof for the cohomological criterion for noetherian regular local rings}

\begin{document}

\date{\today}

\maketitle

\begin{abstract}
Let $(A,\ideal{m}, k=A/\ideal{m})$ be a noetherian local ring. Then it is equivalent
$n = \dim A = \dim_k \ideal{m}/\ideal{m}^2$ and $\Tor^A_i(k,k) = 0$ for all $i \gg 0$.
The article gives a proof with the change of ring spectral sequence in the derived
category form.
\end{abstract}

\section{Theorem and Proof}

\subsection{Acknowledgements}

The author wants to thank Christian Kaiser from the MPIM, Bonn, for proofreading many preliminary
versions of this paper.

\subsection{Auxiliary Results}

We start with some results used later in the proof.

\begin{lemma}

\label{lem:min-free-res-ideal}
Let $(A,\ideal{m})$ be a local ring and $L_\bullet \to M \to 0$ a free $A$--resolution
of a finitely generated $A$--module $M$. Then $(L_\bullet,d_\bullet)$ can be chosen such,
that $d_i(L_i) \subseteq \ideal{m} L_{i-1}$ holds.
\end{lemma}

\small
{\bf Proof.} 
Assume $L_{i-1} \to L_{i-2} \to \cdots \to L_0 \to M$ already constructed and
$0 \to Z \to L_{i-1} \to L_{i-2}$ exact, as well as $L_{i-1}$ with minimal rank. Then it is
$Z \subseteq \ideal{m} L_{i-1}$. 

Write $L_{i-1}$ as $L_{i-1} = A e_1 \oplus \cdots \oplus A e_r$ and let
$z = a_1 e_1 + \cdots + a_r e_r \in Z$. If there is $a_1 \notin \ideal{m}$, therefore
$a_1 \in A^*$, one can assume $a_1 = 1$ without restriction of generality.
So we have the equation $d_{i-1}(e_1) = -a_2 \, d_{i-1}(e_2) - \cdots - a_r \, d_{i-1}(e_r)$
and the summand $A \, e_1$ in $L_{i-1}$ would be superflous and one could chose $L_{i-1}$ with
smaller rank.

At the next step one choses $L_i$ with minimal rank and a surjective map $L_i \to Z \to 0$
and proceeds with the construction inductively.
\normalsize

\begin{definition}
Let $A$ be a commutative ring. Then $\gldh A = \sup \{p \mid \Ext_A^p(M,N) \neq 0\}$,
where $M, N$ are arbitrary $A$--Modules, is called the {\em global dimension} of $A$.
\end{definition}

The following proposition follows from Serre \cite{serre:89}[IV - 35, Corollaire 2]
\begin{proposition}
For a noetherian local ring $(A,\ideal{m}, k=A/\ideal{m})$ we have the equivalence
\begin{enumeratealpha}
\item
It is $\gldh A \leqslant n$.
\item
It is $\Tor^A_{n+1}(k,k) = 0$
\end{enumeratealpha}
\end{proposition}

From Serre \cite{serre:89}[IV - 35, Proposition 21] we cite
\begin{proposition}[Auslander-Buchsbaum-Formula]
Let $(A,\ideal{m})$ be a local noetherian ring with finite global dimension $\gldh A = n$. Furthermore
let $M$ be a finitely generated $A$--Module. Then we have
\begin{equation}
\pd M + \depth M = n = \gldh A
\end{equation}
\end{proposition}
We use this formula several times in the proof without explicitly referring to it.
\smallskip

\subsection{Main Theorem}

\begin{proposition}
Let $(A,\ideal{m})$ be a noetherian local ring of dimension $n$.
Then the following conditions on $A$ are equivalent:
\begin{enumeratealpha}
\item
The homological dimension $\gldh A$ is finite.
In this case it holds automatically that $n = \gldh A = \dim A$.
\item
It is $\dim_k (\ideal{m}/\ideal{m}^2) = n$
\item
$\ideal{m}$ can be generated by $n$ elements.
\item
The tangential cone\index{Tangentialkegel} is $\gr_\ideal{m} A \cong k[X_1,\ldots,X_n]$
\end{enumeratealpha}
where $k=A/\ideal{m}$ is the residue field of $A$.
\end{proposition}
\small
{\bf Proof.}
We will call a ring that fulfills b), c) oder d) {\em geometrically regular}. A
ring that fulfills a) we call {\em cohomologically regular}.

First, because of the Nakayama--Lemma and  because $\ideal{m}$ must be generated by at least
$n$ elements, b) is equivalent to c). Furthermore we have then a surjection
$k[T_1,\ldots,T_n] \twoheadrightarrow \gr_\ideal{m} A$, therefore an isomorphism
$k[T_1,\ldots,T_n]/I \toisom \gr_\ideal{m} A$. As $n = \dim A = \dim \gr_\ideal{m} A$,
it must be $I = 0$. Thereby d) is proved.

Now let d) hold: If $\bar{x}_i$ is the image of $T_i$ in $\ideal{m}/\ideal{m}^2$
and $x_i \in \ideal{m}$ a preimage, then $x_1,\ldots,x_n$ is a regular sequence in $A$.
Therefore $K_\bullet(x_1,\ldots,x_n) \to k \to 0$, the Koszul--complex for the $x_i$,
is a free resolution of length $n$ of $k$. So we have $\Tor_j(k,M) = 0$ for $j > n$ and
$\Tor_n(k,k) = k \neq 0$. This shows $A$ is cohomologically regular with $\gldh A = n$.

\medskip

Now let, for the reverse direction, $(A,\ideal{m})$ be a local noetherian ring with 
$\gldh A < \infty$. We show by induction over $\dim A$, that $A$ is then a geometrically
regular ring too.

\smallskip

\paragraph{Step 1}
First for $\dim A = 0$ the ring $A$ is an artinian ring and $\depth A = 0$. So we have
$\gldh A = \depth A + \pd_A A = 0$ and especially $\dim A = \gldh A = 0$.

It follows that $\pd_A k = 0$ so that $k$ is projective, therefore free, so that $k = A^r$. With
$- \tensor_A k$ it follows $r=1$, so $A=k$ and $\ideal{m} = 0$. This closes the case $\dim A=0$.

\smallskip

\paragraph{Step 2}
Now let $\dim A = n$ and assume the theorem already proven for $\dim A < n$. Choose a prime ideal 
$\ideal{p} \subseteq \ideal{m}$ with $\dim A_\ideal{p} = n-1$ and consider the free
$A$--resolution 
\[
F_\bullet \to A/\ideal{p} \to 0.
\]
It is of finite length. Tensoring by
$A_\ideal{p}$ gives a finite free $A_\ideal{p}$--resolution of 
$k(\ideal{p}) = A_\ideal{p}/\ideal{p} A_\ideal{p}$. Therefore we have 
$\Tor^{A_\ideal{p}}_\nu(k(\ideal{p}),k(\ideal{p}))=0$ for all $\nu \gg 0$. So it is
$\gldh A_\ideal{p} < \infty$ and by the inductive assumption 
$\gldh A_\ideal{p} = \dim A_\ideal{p} = n-1$. The shortest free resolution
\[
F_\bullet \to A/\ideal{p} \to 0
\]
has therefore at least the length $r \geqslant n-1$. Now it is 
$\pd A/\ideal{p} + \depth A/\ideal{p} = \gldh A$. 
Also we have $\depth A/\ideal{p} \geqslant 1$. Therefore $\gldh A \geqslant r + 1 \geqslant n$.

Furthermore $\pd A + \depth A = \gldh A \geqslant n$. As $\pd A = 0$ it is
$\depth A \geqslant n$ and by $\depth A \leqslant \dim A = n$ also $\depth A = \dim A = n$.
So we have $\gldh A = n = \dim A$ and $A$ is a Cohen--Macaulay--ring.

\paragraph{Step 3}
Let $\ideal{p}_i$ be the prime ideals of $\Ass A$, all minimal. Then it holds
$\ideal{m} \not\subset \ideal{m}^2 \cup \bigcup_i \ideal{p}_i$. So a $g \in \ideal{m} - \ideal{m}^2$,
exists, which is not a zero-divisor in $A$ and gives an exact sequence
\[
0 \to A \xrightarrow{\cdot g} A \to A/gA \to 0
\]
We call $A' = A/gA$.

\paragraph{Step 4}

In this step we prove the following lemma, which we formulate explicitly because of its importance
and usefulness:
\begin{lemma}
Let $(A,\ideal{m}, k= A/\ideal{m})$ 
be a local noetherian ring and $g \in \ideal{m} - \ideal{m}^2$ be a non-zero-divisor
of $A$. Furthermore let $A' = A/gA$. Then we have an isomorphism
\begin{equation}
\Tor^{A'}_p(k,k) \oplus \Tor^{A'}_{p-1}(k,k) \cong \Tor^A_p(k,k)
\end{equation}
\end{lemma}
\bigskip 
{\bf Proof.}
We consider in the derived categories of $A$ and $A'$--modules the identity
\begin{equation}
\label{eq:deriv-cat-cr-ss}
k \tensor_{A'}^L (A' \tensor_A^L k) = k \tensor_A^L k
\end{equation}
which is the expression of the change-of-ring spectral sequence in derived category notation.

To understand the following it is just necessary to know, that $M_\bullet \tensor_A^L N_\bullet$
is the total complex of the double complex $P_\bullet \tensor_A N_\bullet$ where 
$P_\bullet \to M_\bullet \to 0$ is a projective $A$-resolution of $M_\bullet$ and a quasi-isomorphism.

To compute $A' \tensor^L_A k$ as an $A'$--complex, choose a free $A$-resolution of $k$, of the form
\[
 \cdots \to F_2 \to F_1 = A^m \to F_0 = A \to k \to 0
\]
and with $d_i(F_i) \subseteq \ideal{m}_A F_{i-1}$,
where $d_1:A^m \xrightarrow{(g, x_2, \ldots, x_m)}  A$ and
$\ideal{m}/\ideal{m}^2 = (g + \ideal{m}^2, x_2 + \ideal{m}^2,\ldots,x_m +\ideal{m}^2)$.
Note that $F_1$ is a minimal rank generator via $d_1$ of the kernel $\ideal{m} = \ker(A=F_0 \to k)$.
So Lemma \ref{lem:min-free-res-ideal} can be applied.

We have then a morphism of complexes of $A'$-modules 
\begin{equation}
\xymatrix@+1.5pc{
 \cdots \ar[r] & F_2 \tensor_A A' \ar[r] \ar[d] \ar[d]^{0} & 
 F_1 \tensor_A A' \ar[r]^{(0,\bar{x}_2,\ldots,\bar{x}_m)} \ar[d]^{(1,0,\ldots,0)}
  & F_0 \tensor_A A' \ar[d]^{1}\\
 \cdots \ar[r] & 0 \ar[r] & k \ar[r]^{0} & k \\
}
\end{equation} 
As the image of $A' \xrightarrow{(1, 0, \ldots, 0)} A'^m = F_1 \tensor_A A'$ is in the kernel of $d_1 \tensor \id_{A'}$ and
not in the image of $d_2 \tensor \id_{A'}$ and
as $\Tor^A_p(A', k) = 0$ for $p > 1$ and $\Tor^A_p(A', k) = k$ for $p \leqslant 1$ by $0 \to A \xrightarrow{\cdot g} A \to A' \to 0$
this gives a quasi-isomorphism of $A'$-complexes
\begin{equation}
A' \tensor_A^L k \toisom \{ \cdots \to 0 \to k \xrightarrow{0} k \to 0\}
\end{equation}

Now chose a resolution
$Q_\bullet \to k \to 0$ of free $A'$--modules with $d_{i}:Q_{i} \to Q_{i-1}$ and 
$d_{i}(Q_i) \subseteq \ideal{m}_{A'} Q_{i-1}$, that is with a matrix representation for all $d_i$
with entries from $\ideal{m}_{A'}$. Calculation in the derived category then gives
\begin{multline}
k \tensor_{A'}^L (A' \tensor_A^L k) = \tot (Q_\bullet \tensor_{A'} \{ 0 \to k \to k \to 0 \}) = \\
= \tot \{ Q_\bullet \tensor_{A'} k \to Q_\bullet \tensor_{A'} k \}
\end{multline}
where in the last complex all derivations are zero. So it is
\begin{multline}
h_p(\tot \{ Q_\bullet \tensor_{A'} k \to Q_\bullet \tensor_{A'} k \}) = 
h_p \left(\tot \left((Q_\bullet \tensor_{A'} k) \oplus (Q_\bullet \tensor_{A'} k)\right)\right) = \\
= h_p(Q_\bullet \tensor_{A'} k) \oplus h_{p-1}(Q_\bullet \tensor_{A'} k) = \Tor^{A'}_p(k,k) \oplus \Tor^{A'}_{p-1}(k,k)
\end{multline}

The right side in \eqref{eq:deriv-cat-cr-ss}
is $P_\bullet \tensor_A k$ where $P_\bullet \to k$ is a free $A$--resolution
of $k$. Comparing cohomology of the two equal (in the derived category) complexes right and left we have
\begin{equation}
\Tor^{A'}_p(k,k) \oplus \Tor^{A'}_{p-1}(k,k) = \Tor^A_p(k,k)
\end{equation}
and the lemma is proven.

\paragraph{Step 5}
Now in our case where $\gldh A < \infty$ we have, that
$\Tor^A_p(k,k)$ vanishes for all $p \gg 0$ and so does $\Tor^{A'}_p(k,k)$ by the lemma of step 4.

So we have $\gldh A' < \infty$ and
by induction, that $A' = A/gA$ is a geometrically regular local ring of dimension $n-1$. 
Therefore $A$ is a geometrically regular local ring of dimension $n$.

\normalsize


\bibliography{mainlitbank}

\end{document}